\documentclass[11pt]{article}
\usepackage[centertags]{amsmath}
\usepackage{amsfonts}
\usepackage{amssymb}
\usepackage{amsthm}
\usepackage{graphicx}

\begin{document}

\title{Cohen-Macaulay graphs arising from digraphs}
\author{GIUSEPPA CARRA' FERRO\\
DANIELA FERRARELLO\\
Department of Mathematics and Computer Science\\
University of Catania\\
Viale Andrea Doria 6, 95125 Catania, Italy\\
e-mail:carra@dmi.unict.it, ferrarello@dmi.unict.it}
\date{}
\maketitle

\newtheorem{defn}{DEFINITION}[section]
\newtheorem{thm}{THEOREM}[section]
\newtheorem{ex}{EXAMPLE}[section]
\newtheorem{rem}{REMARK}[section]
\newtheorem{prop}{PROPOSITION}[section]
\newtheorem{cor}{COROLLARY}[section]

\renewcommand{\qed}{\hfill$\Box$\smallskip}

\begin{abstract}
In this paper we show  a correspondence between directed graphs and
bipartite undirected graphs with a perfect matching, that allows to
study properties of directed graphs through the properties of the
corresponding undirected graphs. In particular it is shown that a
directed graph is transitive iff a corresponding undirected graph is
Cohen-Macaulay.
\end{abstract}


\section{Introduction}
There are several papers on the correspondence between undirected
graphs and polynomial ideals. In particular, there exist many ways
to get properties of an undirected graph from certain binomial and
monomial ideals (\cite{SVV94}, \cite{DS98}, \cite{HO98},
\cite{HO99}, \cite{Ka99}, \cite{BPS01}, \cite{Vi01},
\cite{HH03},\cite{CF06a}). The aim of this paper is to study some
undirected graphs defined as in \cite{CF05c} and associated to a
directed graph $G$ and some relations between their properties.
Furthermore such properties can be verified either by using the
well known algorithms in graph theory or by using the
computational algebra algorithms in polynomial
algebra.\\
After some preliminaries about polynomial toric ideals and graph
theory, in section $4$ we study the undirected graph $H_G$
associated to a digraph $G$. In particular we show a one to one
correspondence between a transitive digraph $G$ and the
Cohen-Macaulay undirected graph $H_G$ associated to it. We also
study the cases when $G$ is a directed acyclic graph (DAG, for
short) by finding the corresponding properties of $H_G$. DAG's play
an important role in graph theory and applications. In fact DAG's
can be used in the study of Bayesian networks, that are directed
graphs with vertices representing variables and edges representing
dependence relations among the variables. Finally digraphs are very
useful  in computational molecular biology (\cite{dJ02}) and in the
minimum cost flow problem in networks, which has many physical
applications (\cite{IsIm00}).

%This last fact allows to have an algorithm for checking when a
%digraph $G$ is Hamiltonian by checking the existence of a perfect
%matching in $K_G$.

\section{Preliminary tools} \label{tools}
\subsection{Toric ideals}

In this section we introduce some properties of the toric
ideals.\\
Let ${\bf N}_{0}$=$\{ 0,1,2,\ldots,n,\ldots \}$ and let
$X_{1},\ldots,X_{n}$ be $n$ variables. Let $K$ be a field of
characteristic zero. Let $A=K[X_{1},\ldots,X_{n}]$ and let
$T_{A}=PP(X_{1},\ldots,\\X_{n})$ be the set of terms of $A$.

\begin{defn}
A term ordering $\sigma$ on $T_{A}$  is a total order such that:\\
(i) $1<_{\sigma} t$ for all $t \in T_{A} \setminus \{1\}$; \\
(ii) $t_{1}<_{\sigma} t_{2}$ implies $t_{1}t'<_{\sigma} t_{2}t'$
for all $t' \in T_{A}$.
\end{defn}

\smallskip
\noindent If $\sigma$ is a term ordering on $T_{A}$ and $f \in A$,
then $M_{\sigma}(f)$ is the monomial $c_{j}t_{j}$ such that
$t_{i}<_{\sigma}t_{j}$ for all $i \neq j$, $i=1,\ldots,r$.
$M_{\sigma}(f)$ is called the {\em leading monomial} of f.

\begin{defn}
Let $\sigma$ be a term ordering on $T_A$ and let $I$ be an ideal
in $A$. The monomial ideal $M_{\sigma}(I)=(M_{\sigma}(f): f \in
I)$ is the {\em initial ideal} of $I$.
\end{defn}

\begin{defn}
{\em (\cite{Bu76})} Let $I$ be an ideal in $A$ and let $\sigma$ be
a term ordering on $T_A$. If $I=(f_{1},\ldots,f_{r})$, then
$\{f_{1},\ldots,f_{r}\}$ is a Gr\"{o}bner basis of $I$ with
respect to $\sigma$ on $T_{A}$ iff
$M_{\sigma}(I)=(M_{\sigma}(f_{1}),\ldots,M_{\sigma}(f_{r}))$.\\
$\{f_{1},\ldots,f_{r}\}$ is a reduced Gr\"{o}bner basis iff
$M_{\sigma}(f_{h})$=$T_{\sigma}(f_{h})$  and $f_h$ is reduced with
respect to $F \setminus f_{h}$  for all $h=1,\ldots,r$.
\end{defn}

\smallskip
\noindent Every ideal $I$ in $A$ has a finite Gr\"{o}bner basis
(\cite{Bu76}). Now we introduce the notion and some properties of
a toric ideal.

\begin{defn}
Let $M$=$(m_{ij})_{i=1,\ldots,m,j=1,\ldots,n}$ be a $(m,n)$-matrix
with $m_{ij}$ in ${\bf Z}$. Let $\pi$ : $K[X_1,\ldots,X_n]
\longrightarrow K[t_1,\ldots,t_m,t_{1}^{-1},\ldots, t_{m}^{-1}]$
be the homomorphism of semigroup algebras defined by
$\pi(X_j)$=$\prod_{i=1,\ldots,m}t_{i}^{m_{ij}}$ for all
$j=1,\ldots,n$. $I_{M}$=$ker(\pi)$ in $A=K[X_1,\ldots,X_n]$ is
called the {\em toric ideal of} $M$.
\end{defn}

\medskip
\noindent
It is well known that $I_M$ is a prime binomial ideal.
Moreover the ideal $(X_{j}-\prod_{i=1,\ldots,m}t_{i}^{m_{ij}},
j=1,\ldots,n)$ is  in
$K[X_1,\ldots,X_n,t_1,\ldots,t_m,t_{1}^{-1},\ldots, t_{m}^{-1}]$
and $I_{M}=(X_{j}-\prod_{i=1,\ldots,m}t_{i}^{m_{ij}},
j=1,\ldots,n) \cap K[X_1,\ldots,X_n]$ (\cite{St95}).

\begin{thm}
{\em \cite{St95}}. $I_M$ is generated by binomials of the type
$X^{u^{+}}-X^{u^{-}}$, where $u^{+}, u^{-} \in {\bf Z}^{n}$ are
non negative with disjoint support.
\end{thm}

\subsection{Graphs and digraphs}

In this paper $G$=$(V(G),E(G))$ will be a finite graph with
$V(G)=\{v_1,\ldots,v_n\}$ and $E(G)=\{e_1,\ldots,e_m\}$.
Furthermore  $[v_{i},v_{j}]$ will denote the directed edge from
$v_{i}$ to
$v_{j}$, while $\{v_i,v_j\}$ will denote the undirected edge between $v_{i}$ and $v_{j}$.\\
The {\em underlying graph} $G_u$ of a directed graph $G$ is the
undirected graph with $V(G_u)=V(G)$ and the same undirected edges
of $G$. All graphs in this paper will be simple, i.e. without
multiple edges. A directed graph $G$ will be shortly called a {\em
digraph}.

%\begin{defn}
%Let $G$ be a graph (respectively a digraph). A {\em walk} of length
%$n$ from a vertex $v_i$ to a vertex $v_j$ in $G$ is a sequence of
%vertices $v_i$=$v_{i(1)}$, \ldots, $v_j$=$v_{i(n+1)}$, such that
%$\{v_i(h),v_i(h+1)\}$ (respectively either $[v_{i(h)}v_{i(h+1)}]$ or
%$[v_{i(h+1)}v_{i(h)}]$) is in $E(G)$) for all $h=1,\ldots,n$. \\
%A walk in $G$ is called {\em simple} if there are no repeated edges.
%A walk in $G$ is called {\em elementary} if there are no repeated
%vertices. A walk in a digraph $G$ is called a {\em directed walk} if
%$[v_{i(h)}v_{i(h+1)}]$ is in $E(G)$ for all $h$. If
%$v_{i(1)}=v_{i(n+1)}$, and the walk is elementary, then it is a {\em
%cycle} (respectively a {\em directed cycle}).
%\end{defn}

%\smallskip
%\noindent Here we introduce other definitions in graph theory and
%some corresponding properties.
\begin{defn}
An undirected graph $G$ is {\em bipartite} if its vertices can be
divided in two sets, such that no edge connects vertices in the
same set. Equivalently $G$ is {\em bipartite} iff all cycles in
$G$ are even.  $G$ is {\em acyclic} if it has no cycle.
\end{defn}

\begin{defn}
A digraph $G$ is called a {\em directed acyclic graph}, {\em DAG}
for short, when there are no directed cycles in $G$.
\end{defn}

\begin{defn}
A digraph $G=(V,E)$ is called {\em transitive} if for all $u$,
$v$, $w$ $\in V$, such that $[u,v] \in E$ and $[v,w] \in E$ we
have $[u,w] \in E$.
\end{defn}

\smallskip
\noindent Equivalently $G$ is a transitive graph if there exist a
directed edge from $u$ to $v$, whenever there is a directed walk
from $u$ to $v$.

\begin{rem}
A simple directed transitive graph is a DAG. In fact every possible cycle has
to be undirected.
\end{rem}

%\begin{defn}
%Let $G$ be a digraph. A vertex $v_i$ in $V(G)$ is called a {\em source} if no adjacent
%edge is directed into $v_i$. A vertex $v_i$ in $G$ is called a {\em sink} if every adjacent edge is directed into %%@
%$v_i$.
%\end{defn}

\begin{defn}
A {\em vertex cover} $V'$ of an undirected graph $G$ is a subset of
vertices of $G$, such that at least one of the vertices of
every edge in $G$ is in $V'$. A vertex cover $V'$ is said
to be {\em minimal} if no subsets of $V'$ is a vertex cover.
Finally a graph is said {\em unmixed} if
all minimal vertex covers have the same cardinality.
\end{defn}

%\smallskip
%\noindent A generalization of the concept of a minimal vertex
%cover for a digraph is the following one.

%\begin{defn}
%Let $G=(V,E)$ be a digraph. A {\em source cover of} $G$ is a
%vertex cover $V'$ of $G$, such that every edge in $G$ leaves every
%vertex in $V'$. A source cover $V'$ of $G$ is called {\em minimal}
%if no subset of $V'$ is a source cover of $G$. A {\em sink cover
%of} $G$ is a vertex cover $V'$ of $G$, such that every edge in $G$
%does not leave any vertex in $V'$. A sink cover $V'$ of $G$ is
%called {\em minimal} if no subset of $V'$ is a sink cover of $G$.
%\end{defn}

\section{Ideals arising from graphs and digraphs}
Here we will introduce the binomial and monomial ideals, that are
associated to an undirected   graph (respectively a digraph) as in
(\cite{CF06a}) and (\cite{CF05c}). All ideals in the paper are in a
polynomial ring with coefficients in a field $K$ of characteristic
zero.

\begin{defn}
Let $G$ be an  undirected graph. The binomial {\em extended edge
ideal} of $G$ is the ideal $I(G,E(G))$=$($ $e_h-v_{i}v_{j}$:
$e_h=\{v_{i},v_{j}\} \in E(G)$, $i=1,\ldots,n$ $)$ in
$K[e_1,\ldots,e_m,v_1,\ldots,v_n,z_{1},\ldots,z_{n}]$. The ideal
$I(G)_{E(G)}=I(G,E(G))\cap K[e_1,\ldots,e_m]$ is the {\em binomial
edge ideal of} $G$.
\end{defn}

\begin{rem}
$I(E)_G$ is the toric ideal of the matrix $IM(G)^{t}$, that is the
incidence matrix $IM(G)=(a_{ih})_{i=1,\ldots,n,h=1,\ldots,m}$ of $G$
defined by $a_{ih}=1$ if $v_{i} \in e_{h}$ and $a_{ih}=0$ if $v_{i}
\notin e_{h}$ for every $v_{i} \in V(E)$ and $e_{h} \in E(G)$. Its
definition can be also found in $(\cite{HO98})$ and $(\cite{Ka99})$.
\end{rem}

\smallskip
\noindent This notion can be extended to digraphs.

\begin{defn}
Let $G$ be a digraph. The binomial {\em extended diedge ideal} of
$G$ is the ideal $I(G,E(G))$=$($ $e_h-z_{i}v_{j}, z_{i}v_{i}-1$:
$e_h=[v_{i},v_{j}] \in E(G)$, $i=1,\ldots,n$ $)$ in
$K[e_1,\ldots,e_m,v_1,\ldots,v_n,z_{1},\ldots,z_{n}]$. The ideal
$I(G)_{E(G)}=I(G,E(G))\cap K[e_1,\ldots,e_m]$ is the {\em binomial
diedge ideal of} $G$.
\end{defn}

\begin{rem}
$I(G)_{E(G)}$ is the toric ideal of the matrix $IM(G)^{t}$, that
is the transpose of the incidence matrix
$IM(G)=(a_{ih})_{i=1,\ldots,n,h=1,\ldots,m}$ of $G$ defined by
$a_{ih}=-1$ if $e_{h}$ leaves $v_{i}$, $a_{ih}=1$ if $e_{h}$
arrives to  $v_{i}$ and $a_{ih}=0$ if $v_{i} \notin e_{h}$ for
every $v_{i} \in V(E)$ and $e_{h} \in E(G)$.
\end{rem}

\medskip
\noindent Given an undirected graph $G$ and  an even closed walk \ \
$C=(e_{i_{1}}=\{v_{i_{1}},v_{i_{2}}\}$,
$e_{i_{2}}=\{v_{i_{2}}$,$v_{i_{3}}\}$, $\ldots$,
$e_{i_{2q-1}}$=$\{v_{i_{2q-1}}$,$v_{i_{2q}}\}$,
$e_{i_{2q}}$=$\{v_{i_{2q}}$,$v_{i_{1}}\})$ of $G$, let \\
$f_{C}$=$\prod_{k=1,\ldots,q}e_{i_{2k-1}}-\prod_{k=1,\ldots,q}e_{i_{2k}}$
be the corresponding binomial in $I(G)_{E(G)}$.\\
Given a digraph $G$ and a closed   walk
$C=(e_{i_{1}}=[v_{i_{1}},v_{i_{2}}],
e_{i_{2}}=[v_{i_{2}},v_{i_{3}}], \ldots,\\
e_{i_{q-1}}=[v_{i_{q-1}},v_{i_{q}}],
e_{i_{q}}=[v_{i_{q}},v_{i_{1}}])$ of $G$, let $f_{C}=\prod_{i \in
I}e_{i}-\prod_{j \in J}e_{j}$, be the corresponding binomial in
$I(G)_{E(G)}$ with $I,J \subset \{1,\ldots,n\}$ and $|I|+|J|=q$.

\smallskip
\noindent If $G$ is an undirected graph, then the toric ideal
$I(G)_{E(G)}$ is generated by all binomials $f_C$, where $C$ is an
even closed walk of $G$ (\cite{Vi95} and \cite{HO99}). If $G$ is a
digraph, then toric ideal $I(G)_{E(G)}$ is generated by all
binomials $f_{C}$, where $C$ is a cycle of $G$ (\cite{CF05c}).

\smallskip
\noindent Now we introduce some binomial ideals associated to the
vertices of a graph and a digraph.

\begin{defn}
{\em (\cite{CF06a})}\\Let $G$=$(V(G),E(G))$ be a graph. The ideal
$I(G,V(G))$=$($$v_i-\prod e_{h}$: $v_i$ belongs to the edge
$e_{h}$$)$ in $K[v_1 ,\ldots, v_n, e_1, \ldots, e_m]$ is the
{\em binomial extended vertex ideal of} $G$.\\
$I(G)_{V(G)}$=$I(G,V(G)) \cap K[v_1, \ldots, v_n$ , $z_1, \dots,
z_n]$ is the {\em vertex ideal of} $G$.
\end{defn}

\begin{defn}
{\em (\cite{CF05c})}\\ Let $G$=$(V(G),E(G))$ be a finite digraph.
The ideal $I(G,V(G))$=($v_{i}-\prod e_{h}$: $e_h$ arrives in $v_i$,
$z_{i}-\prod e_{h}$: $e_h$ leaves from $v_i$) in $K[v_1,\ldots,v_n,
z_1, \ldots, z_n, e_1,\ldots,e_m]$ is the binomial {\em extended
divertex ideal of} $G$.\\ $I(G)_{V(G)}$=$I(G,V(G))\cap K[v_1,\ldots,
v_n$, $z_1, \dots, z_n]$ is the {\em divertex ideal of} $G$.
\end{defn}

\section{The graph $H_G$}
\smallskip
\noindent The following undirected graph  associated to a digraph
has some nice properties.
\begin{defn}
Let $G$ be a digraph.  Let $H_G$ be the undirected graph with
$V(H_G)$=$V(G) \cup \{z_{1},\ldots,z_{n}\}$ =
$\{v_{1},\ldots,v_{n},z_{1},\ldots,z_{n}\}$
 and
 $E(H_G)$=\\\{$e_h=\{z_{i},v_{j}\}$: $\{v_{i},v_{j}\}
\in E(G)$, $h=1\ldots,m$\} $\cup$ \{$f_{i}=\{z_{i},v_{i}\}$:
$i=1,\ldots,n\}$. Let $R=K[v_{1},z_{1},f_{1}, \ldots,
v_{n},z_{n},f_{n}, e_{1}, \ldots, e_{m}]$ and let  $\pi: R
\rightarrow R/(f_{1}-1,\ldots f_{n}-1)$ be the canonical ring
homomorphism defined by $\pi(v_{i})=v_{i}$, $\pi(z_{i})=z_{i}$,
$\pi(f_{i})=1$, $\pi(e_{j})=e_{j}$ for all $i=1,\ldots,n$ and
$j=1,\ldots,m$. $\pi(I(H_G,E(H_G))=I(G,E(G))$.
\end{defn}

It is proved in \cite{CF05c} that an undirected graph $G'$ is
bipartite with a perfect matching  if and only if there exists a
digraph $G$ without loops, such that $G'=H_G$.

\begin{rem}
$G$ is not unique. If we consider the graph $G^S$, {\em edge
symmetric} of $G$, i.e. the graph with $V(G^S)=V(G)$ and $[v_i,v_j]
\in E(G^S)$ if and only if $[v_j,v_i] \in E(G)$, then $H_{G^S}=G'$.
\end{rem}

\smallskip
\noindent
 Now we study some properties of $H_G$ that are in correspondence with
properties of $G$ and depend on the cycles in both of them.

\bigskip
\noindent We start with the case  when $G$ has no directed cycles
and the relation with undirected Cohen-Macaulay graphs.

\smallskip
\noindent
To every undirected graph $G$ with vertex set
$V(G)=\{v_1, \ldots, v_n\}$ and edge set $E(G)=\{e_1,\ldots,e_m\}$
it is possible to associate a monomial ideal $I(G)$, that is
generated by all square free monomials $v_i v_j$, such that
$\{v_i,v_j\}=e_h$ is an edge of $G$. This ideal is usually
called the monomial {\em edge ideal}.

\begin{defn}
An undirected graph $G$ is  {\em Cohen-Macaulay} (CM for short)
with respect to the field $K$, if the quotient ring $K[v_1,
\ldots, v_n]/I(G)$ is a Cohen-Macaulay ring.
\end{defn}

\smallskip
\noindent
It is known that every CM graph is unmixed.\\
CM graphs are studied in several works (\cite{Vi90}, \cite{Vi01},
\cite{HH03}, \cite{HHX04}, \cite{Fa04}, \cite{Fa05}, \cite{FV05}).
Actually there is no general decision procedure for CM graph just by
using only combinatoric tools. There is a combinatoric
characterization of bipartite CM graphs in \cite{HH03} and a
decision procedure for bipartite CM graphs can be found in
\cite{CF04b}.

\smallskip
\noindent
Here we will investigate the case when $H_G$ is Cohen-Macaulay.\\
We will show that the condition on $G$ for $H_G$ to be CM is that
$G$ is a transitive graph, while a weaker condition than
Cohen-Macaulayness for $H_G$ is equivalent for $G$ to be a DAG. In
order to show these results we need the following useful theorem
by Herzog and Hibi, that characterizes bipartite CM graphs.

\begin{thm} \label{bipCM} {\em (\cite{HH03})}\\
Let $G'$ be a  simple bipartite graph without loops on the vertex
set $W\bigcup W^\prime$, with $W=\{x_1, \ldots, x_n\}$, $W^\prime=\{y_1,
\ldots, y_n\}$ such that
\begin{description}
    \item[(a)]$\{x_i,y_i\}$ is an edge for all $1\leqslant i \leqslant n$;
    \item[(b)]if $\{x_i,y_j\}$ is an edge then $i \leqslant j$;
\end{description}
then $G'$ is CM iff
\begin{description}
    \item[(c)] whenever $\{x_i,y_j\}$ and $\{x_j,y_k\}$ are edges, then $\{x_i,y_k\}$ is
an edge.
\end{description}
\end{thm}

\begin{rem} \label{a condition}
Since $G'$ is  bipartite with the bipartition sets $W=\{v_1,
\ldots, v_n\}$ and $W'=\{z_1, \ldots, z_n\}$, then    the condition
\textbf{(a)} is equivalent to the existence of the perfect matching
$M$=\{$\{v_i,z_i\}$: $i=1,\ldots,n$\} in $G'$.\\
If $G'$=$H_G$ for some digraph $G$, then the condition \textbf{(a)} is always true.
\end{rem}

\smallskip
\noindent We have the following results on DAG's and transitive
graphs.

\begin{prop}
Let $G$ be a digraph with $|V(G)|=n$. If $G$ is a DAG, then it has a
{\em topological sort}, i.e. a total ordering of the vertices
$\{v_{i(1)},\ldots,v_{i(n}\}$, such that each edge in $G$ is of the
kind $[v_{i(1)},v_{i(j)}]$ with $i<j$.
\end{prop}

\smallskip
\noindent
 In general, this ordering is not unique (\cite{Sk90}).

\begin{thm} \label{dag}
Let $G=(V(G),E(G))$ be a simple digraph. Then condition
\textbf{(b)}  for $H_G$ is equivalent to say that $G$ is a DAG.
\end{thm}
\textbf{Proof}. Suppose that  $H_G$ satisfies \textbf{(b)} and $G$
is not a DAG. Then $G$ contains at least a directed cycle
$C$=$\{[v_{i1},v_{i2}]$, $[v_{i2},v_{i3}]$, $\ldots$,
$[v_{ik},v_{i1}]\}$. So, the set $L$=\{$\{v_{i1},z_{i1}\}$,
$\{z_{i1},v_{i2}\}$, $\{v_{i2},z_{i2}\}$, $\{z_{i2},v_{i3}\}$,
$\dots$, $\{v_{ik},z_{ik}\}$, $\{z_{ik},v_{i1}\}$\}  is a cycle in
$H_G$ and it is impossible to relabel the vertices in order to
satisfy  condition \textbf{(b)}. Now let $G$ be a DAG with
$|V(G)|=n$. By proposition 4.1 we can relabel the vertices $\{v_1,
\ldots, v_n\}$ of $G$ in such a way as $[v_i,v_j]$ $\in E(G)$ if
and only if $i<j$, for all $i,j \in V(G)$. By using this ordering
$\{z_i,v_j\}$ $\in E(H_G)$ if and only if $i\leq j$. So we have
condition \textbf{(b)}. \qed

\begin{thm}\label{hgCM}
Let $G=(V(G),E(G))$ be a simple digraph. The following facts are
equivalent:
\begin{enumerate}
\item $G$ is transitive
\item $H_G$ is Cohen-Macaulay
\end{enumerate}
\end{thm}
\textbf{Proof}. Let $G$ be a directed simple graph. First suppose
that $G$ is transitive. By remark \ref{a condition}, condition
\textbf{(a)} of theorem \ref{bipCM} is satisfied. Condition
\textbf{(b)} is true, by theorem \ref{dag}, since every simple
transitive graph is a DAG, as already observed. Condition
\textbf{(c)} follows by definition of transitivity. Conversely
suppose that $H_G$ is CM. Condition \textbf{(c)} implies that $G$
is transitive.\qed

The construction of $H_G$ and the previous theorem are used in
\cite{VV07}.

\begin{rem}
In \cite{CF06b} it is shown a decision procedure for checking when a
bipartite undirected graph $G$ is CM. By using the same procedure we
have a decision procedure for checking when a digraph is transitive,
without using the other well known algorithms in graph theory
(\cite{Sk03} and \cite{Tr99}). Conversely we can use the well known
algorithms in graph theory for checking the transitivity of a
digraph in order to check the Cohen-Macaulay property of a bipartite
undirected graph.
\end{rem}

\section{The graph $K_G$}

\smallskip
\noindent Another undirected graph associated to a digraph is now
introduced.
\begin{defn}
Let $G$ be a digraph. Let $K_G$ be the undirected graph with
$V(K_G)$=$V(G) \cup
\{z_{1},\ldots,z_{n}\}$=$\{v_{1},\ldots,v_{n},z_{1},\ldots,z_{n}\}$
 and $E(K_G)$=\{$e=\{z_{i},v_{j}\}$: $[v_{i},v_{j}]
\in E(G)$\}. $K_G$ is called the {\em sink-source undirected graph}
associated to $G$.
\end{defn}

\smallskip
\noindent $K_G$ is a subgraph of $H_G$ by its own definition.

\begin{rem} \label{v ideal KG}
If $G$ is a digraph, and $G^*=K_G \setminus L$, where $L$ is the set
of isolated vertices, then the extended divertex ideal of $K_G$ is
equal to $I(G^*,V(G^*))$.
\end{rem}

\subsection{The  Cohen-Macaulay property for $K_G$}

Since we know the characterization  of the Cohen-Macaulay property
for $H_G$  trough the corresponding properties of the digraph $G$ as
in theorem \ref{hgCM} and $K_G$ is a subgraph of $H_G$, it is
natural to ask whether this property is preserved in $K_G$.

\medskip
\noindent The answer in general is negative. In fact when $H_G$ is
C-M graph there are several cases when $K_G$ is C-M and when $K_G$
is not C-M as in the examples below.

\begin{rem}
A necessary condition on $G$ for the Cohen-Macaulay property on
$K_G$ is \\
\textbf{(1)} The number of vertices that are pure sources is equal
to the number of vertices that are pure sinks. \\
In fact every
source and every sink in $G$ determines an isolated vertex in $K_G$.
Now the cardinalities of the two bipartition sets in every connected
component in $K_G$ have to be equal, so it is necessary to have the
same number of isolated vertices among z's and v's.
\end{rem}

\begin{ex}
Let $G_1$ be the digraph with $V(G_1)=\{v_1, \ldots, v_4\}$ and
$E(G_1)=\{[v_1,v_2],[v_1,v_3],[v_1,v_4],[v_2,v_3],[v_2,v_4]\}$.
$G_1$ is transitive, so $H_{G_1}$ is C-M, by \ref{hgCM}. $K_{G_1}$
is given by $V(K_{G_1})=\{z_1, \ldots, z_4, v_1, \ldots, v_4\}$
and
$E(K_{G_1})=\{\{z_1,v_2\},\{z_1,v_3\},\{z_1,v_4\},\{z_2,v_3\},\{z_2,v_4\}\}$,
$K_{G_1}$ is not C-M. In fact it has four connected components:
three isolated vertices ($v_1$, $z_3$ and $z_4$), that are
trivially Cohen-Macaulay, and a bipartite graph with $5$ vertices,
that is not C-M by \ref{bipCM}.
\end{ex}

\begin{ex}
Let $G_2$ be the digraph with $V(G_2)=\{v_1, \ldots, v_4\}$ and
$E(G_2)=\{[v_1,v_2],[v_1,v_4],[v_3,v_2],[v_3,v_4]\}$. $G_2$ is
transitive, so $H_{G_2}$ is C-M, by \ref{hgCM}. $K_{G_2}$ is given
by $V(K_{G_2})=\{z_1, \ldots, z_4, v_1, \ldots, v_4\}$ and
$E(K_{G_2})=\{\{z_1,v_2\},\{z_1,v_4\},\{z_3,v_2\},\{z_3,v_4\}\}$,
$K_{G_2}$ is not C-M. In fact it has four connected components: four
isolated vertices ($v_1$, $v_3$, $z_2$ and $z_4$), that are
trivially Cohen-Macaulay, and a bipartite graph with $4$ vertices,
that is not C-M by \ref{bipCM}. $G_2$ is transitive and $K_{G_2}$
verifies the condition \textbf{(1)} expressed in the previous
remark, but once again it is not C-M.
\end{ex}

\begin{ex}
Let $G_3$ be the digraph with $V(G_3)=\{v_1, \ldots, v_5\}$ and
$E(G_3)=\{[v_1,v_2],[v_1,v_3],[v_1,v_5],[v_2,v_3],[v_2,v_5],[v_4,v_5]\}$.
$G_3$ is transitive, and so $H_{G_3}$ is Cohen-Macaulay. $K_{G_3}$
has four isolated vertices (two pure sources and two pure sinks) and
a connected component whose vertices are $\{z_1, z_2, z_4, v_2, v_3,
v_5\}$ and whose edges are $\{\{z_1,v_2\}$, $\{z_1,v_3\}$,
$\{z_1,v_5\}$, $\{z_2,v_3\}$, $\{z_2,v_5\}$, $\{z_4,v_5\}\}$, that
is C-M by \ref{bipCM}.
\end{ex}

\medskip
\noindent
 Here we show some necessary conditions, that allow to
prove the C-M property for the undirected graph $K_G$ associated
to a digraph $G$.

\begin{prop}
Let $G$ be a transitive digraph with vertex set
$V(G)=\{v_1,\ldots,v_n\}$, satisfying \textbf{(1)} and let $K_G$ be
the associated source-sink undirected graph. Let us suppose $K_G$
has a perfect
matching $M$ and condition \textbf{(b)} in \ref{bipCM} is not satisfied.\\
Then there exists a perfect matching $M_1$ for $K_G$, with $M_1 \neq
M$.
\end{prop}

\textbf{Proof}. The graph $K_G$ has a perfect matching, by
hypothesis, and it is bipartite, by definition, so thanks to \ref{HG
of a digraph} there exists a digraph $G'$, such that $K_G=H_{G'}$.
Moreover condition \textbf{(b)} is not verified for $K_G$, so by
virtue of \ref{dag} $G'$ is not a dag. Then there exists a directed
cycle in $G'$ of length $k$, that determines a cycle $C$ in $K_G$ of
length $2k$; let $C=\{e_{i1}, e_{i2}, \ldots, e_{i 2k}\}$. By
construction (see \cite{CF05c} teorema sui cicli (thm 5)), $C$ is
such that $k$ disjoint edges of $C$ are contained in the matching
$M$. Then the set of edges of $C$ is $N_1 \cup N_2$, with $N_1
\subseteq M$.

Let $M_1=N_2 \cup (M \backslash N_1)$. $M_1$ is a perfect matching
for $K_G$. In fact every edge in  $N_2$ is not in $M$, because it
shares its vertices with edges in $N_1$, so in $M$, that is a
matching. So, the edges in $M_1$ are disjoint and $M_1$ is a
matching. Finally, since $N_2 \cap M =\emptyset$, a fortiori $N_2
\cap (M \backslash N_1)= \emptyset$, so $|M_1| = |N_2|+|M-N_1| =
k+(n-k) = n$; this suffices to state that $M_1$ is perfect.

 \qed

\begin{cor}
Let $G$ be a transitive digraph with vertex set
$V(G)=\{v_1,\ldots,v_n\}$ and let $K_G$ be the associated
source-sink undirected graph. If $K_G$ has a unique perfect
matching, then the conditions \textbf{(a)} and \textbf{(b)}
described in \ref{bipCM}, in order to be Cohen-Macaulay are
satisfied.
\end{cor}

\begin{prop}
Let $G$ be a transitive digraph with vertex set
$V(G)=\{v_1,\ldots,v_n\}$ and let $K_G$ be the associated
source-sink undirected graph.  $K_G$ is CM if and only if the
following conditions are satisfied.
\begin{description}
\item[(1)] $s$=$sink(G)$=$|\{$ $v \in V(G)$ : v is a sink in $G$
$\}|$=$source(G)$=$|\{$ $v \in V(G)$ : v is a source in $G$ $\}|$;
\item[(2)] $|E(G)| \leq \frac{(n-s+1) (n-s)}{2}$;
 \item[(3)]
$K_G$ has a connected component $G'$ with $|E(G')|$=$|K_G|$ and $G'$
has a unique perfect matching $M$;
\item[(4)] if
$[v_i,v_j],[v_h,v_j], [v_h,v_k] \in E(G)$, with $\{z_h, v_j\} \in
M$, then $[v_i,v_k] \in E(G)$ for all $i,j,h,k=1,\ldots,n$, $i \neq
j \neq h \neq k$.
\end{description}
\end{prop}
\textbf{Proof}. If $K_G$ is CM, then (1), (2), (3) and (4) are
satisfied. In fact since $K_G$ is a subgraph of the bipartite
graph $H_G$, then it is bipartite by its own definition. In fact
$K_G$ cannot have odd cycles. $K_G$ CM implies that it has the
connected components $\{$ $v \in V(G)$ : v is a sink in $G$ $\}$,
$\{$ $v \in V(G)$ : v is a source in $G$ $\}$ and the connected
component $G'$ with $V(G')$=$\{$ $z_{\sigma(j)}, v_{\tau(j)}$:
$j=1,\ldots,n-s$, $v_{\sigma(j)}$ source in $G$ and $v_{\tau(j)}$
sink in $G$ for all $j=n-s+1,\ldots,n$, $\sigma, \tau \in S_n$
$\}$ by (a) of \ref{bipCM}. Since $G'$ has to be bipartite and CM,
then $|V(G')|$ is even and then (1) is satisfied. (2) follows from
the definition of $K_G$. In fact $|E(G)|$=$|E(K_G)|$ and (2) must
be satisfied, because a bipartite CM graph with $2(n-s)$ vertices
has at most $\frac{(n-s+1) (n-s)}{2}$ edges by (a), (b) and (c) of
\ref{bipCM}. (3) is satisfied because a bipartite CM graph has a
perfect matching by (a) of \ref{bipCM} and such matching is unique
because it is a comparability graph by (b) and (c) of \ref{bipCM}.
Finally (4) is equivalent to the property (c). In fact by
\ref{bipCM} we can suppose that $E(K_G)$=$\{$ $[z_{\sigma(i)},
v_{\tau(j)}]$: $i,j=1,\ldots,n-s$ and $i \leq j$ $\}$ and (4) is
equivalent to the property (c) with $\sigma(i)=i$,
$\sigma(j)=h$, $\tau(j)=j$, $\tau(h)=k$.\\
Conversely suppose that $G$ satisfies conditions (1), (2), (3) and
(4). Condition (1) implies that $K_G$ is a subgraph of $H_G$ with
connected components $\{$ $v \in V(G)$ : v is a sink in $G$ $\}$,
$\{$ $v \in V(G)$ : v is a source in $G$ $\}$ and the connected
component $G'$ with set of vertices $V(G')$=$\{$ $z_{\sigma(j)},
v_{\tau(j)}$: $j=1,\ldots,n-s$, $v_{\sigma(j)}$ source in $G$ and
$v_{\tau(j)}$ sink in $G$ for all $j=n-s+1,\ldots,n$, $\sigma, \tau
\in S_n$ $\}$. $E(K_G)$=$E(G')$ by definition of $K_G$ and $G'$ is a
bipartite subgraph of $H_G$ with $2(n-s)$ vertices and bipartition
set $(Z,V)$, with $Z$=\{ $z_{\sigma(j)}$: $j=n-s+1,\ldots,n$ \} and
$V$=\{ $v_{\tau(j)}$: $j=n-s+1,\ldots,n$ \}. Since every vertex in
$G'$ is in some edge in $K_G$ by definition of $K_G$ and by (1),
then $G'$ bipartite implies by Hall's theorem that it has a perfect
matching and the edge set $E(G')$ is union of disjoint matchings of
$G'$. Furthermore such perfect matching is unique by (3). By (2)
$|E(G)|$=$|E(K_G)|$=$|E(G')|$ satisfies the necessary condition on
the maximal number of edges in a CM bipartite graph. By (3) we can
suppose that $M$=\{$[z_{\sigma(j)},v_{\tau(j)}]$: $j=n-s+1,\ldots,n$
\} is the perfect matching of $G'$ and
$M'$=\{$[z_{\sigma(j)},v_{\tau(j)}]$: $j=n-s+2,\ldots,n$ \} is not a
perfect matching of $G'$. Since $M$ is unique, then it follows that
there is no edge $[z_{\sigma(j)},v_{\tau(n-s+1)}]$,
$j=n-s+2,\ldots,n$ in $E(G')$. \qed

\begin{rem}\textbf{Directed paths in $G$}\\
Let $G$ be a transitive digraph on the vertex set $V=\{v_1, \ldots,
v_n\}$, with $s$ number of sources and sinks and let $M$ be the
unique perfect matching in $K_G$. Then it is possible to say
something on some directed paths in $G$, coming from $M$.\\
Since $G$ is transitive, it is a DAG and so it is possible to
reorder the vertices in such a way by drawing them left to right,
all the edges are directed from left to right (topological sort). In
general there are several of these orderings on $V$, but there is at
least one of them, namely $O=\{v_{i1}, \ldots, v_{in}\}$, such that
the edges coming from $M$ are directed from $v_{ij}$ to $v_{i j+1}$;
in fact it is sufficient to choose the vertices in $O$ by using the
rule that $v_i$ is on the left on $v_j$ if the edge $\{z_i, v_j\}$
is in the matching $M$. This rule is consistent with our graph,
because it can not happen to find a vertex on the left of itself,
otherwise the graph $G$ should have a directed cycle, against the
hypothesis of transitivity. So, there exist in $G$ a set $P$ of $s$
directed paths, drawn from left to right. In particular every path
in $P$ starts in one of the $s$ source and finishes in one of the
$s$ sinks. Moreover, the set $P$ mentioned above touches all the
vertices in $V$, so it is a particular edge cover for $G$.
\end{rem}

\begin{ex}
Let $G_4$ be the digraph on the vertex set $\{v_1, \ldots, v_5\}$
and edges $\{[v_2,v_1],[v_1,v_5],[v_2,v_5],[v_2,v_3],[v_4,v_3]\}$.\\
$G_4$ is transitive, and it has two sources ($v_2$ and $v_4$) and
two sinks ($v_3$ and $v_5$). The associated sink-source graph
$K_{G_4}$ has four isolated vertices ($v_2, v_4, z_3, z_5$) and a
unique perfect matching $M=\{\{z_2,v_1\},\{z_1,v_5\},\{z_4,v_3\}\}$.
Among the possible topological sorts over the vertices, we can
choose our ordering in such a way $v_2<v_1, v_1<v_5$ and $v_4<v_3$.
By taking $O=\{v_2<v_1<v_5<v_4<v_3\}$, it is possible to draw the
vertices in such a way all the edges are directed from left to right
and the directed paths $p_1=\{v_2,v_1,v_5\}$ and $p_2=\{v_4,v_3\}$
can be drawn in such a way every edge goes from a vertex to the
following one.
\end{ex}

\bigskip
\noindent %By using the packages {\em networks} and {\em Groebner} of
%Maple 10, we implemented procedures that, given a digraph $G$ allow
%to decide if $G$ has a sink-source vertex cover, by finding
%the set of minimal vertex covers of $K_G$. Moreover if $K_G$ has
%not a perfect matching it is possible to deduce that $G$ is not
%Hamiltonian.
The problem of finding a perfect matching in a undirected graph can
be solved by finding the complement of the minimal edge covers of
$K_G$, that coincide with the minimal vertex cover of the edge
hypergraph associated to $K_G$. Another approach for the perfect
matching is the Hungarian method in the package {\em Combinatorica}
of Mathematica 4.2 .

\end{document}